\documentclass[11pt]{article}

\usepackage{latexsym, amssymb}

\newcommand{\be}{\begin{equation}}
\newcommand{\ee}{\end{equation}}
\newcommand{\beq}{\begin{eqnarray}}
\newcommand{\eeq}{\end{eqnarray}}

\newtheorem{thm}{Theorem}
\newtheorem{prop}{Proposition}
\newtheorem{cor}{Corollary}

\newtheorem{df}{Definition}

\def\stop{\hfill$\Box$}

\def\S{\Sigma}
\def\pf{\noindent {\em{Proof}}: }
\def\C{\mathcal{C}}
\def\R{\mathbb{R}}

\def\Sp{\Sigma^\prime}
\def\G{\mathcal{G}}

\def\mpz{\bar{m}_0}
\def\Abz{\bar{A}_0}
\def\rbz{\bar{r}_0}
\def\ep{\epsilon}

\def\vh{\vspace{.3cm}}

\begin{document}

\title{On The Capacity of Surfaces in Manifolds with Nonnegative Scalar Curvature}

\author{ 
Hubert Bray\thanks{%
Department of Mathematics,
Duke University, Durham, NC 27708, USA. 
E-mail: bray@math.duke.edu}
\ and
Pengzi Miao\thanks{%
School of Mathematical Sciences, Monash University, Victoria, 3800,
Australia.  E-mail: Pengzi.Miao@sci.monash.edu.au}
}

\date{}

\maketitle 

\begin{abstract}
Given a surface in an asymptotically flat $3$-manifold with nonnegative scalar curvature, 
we derive an upper bound for the capacity of the surface in terms of the area of the 
surface and the Willmore functional of the surface. The capacity of a surface is 
defined to be the energy of the harmonic function which equals $0$ on the surface 
and goes to $1$ at $\infty$. Even in the special case of $\R^3$, this is a new estimate. 
More generally, equality holds precisely for a spherically symmetric sphere in 
a spatial Schwarzschild $3$-manifold. As applications, we obtain inequalities relating the 
capacity of the surface to the Hawking mass of the surface and the total mass of 
the asymptotically flat manifold. 

\end{abstract}

\section{Introduction}

The research in this paper was partly motivated by the following theorem.

\vh

\noindent {\bf Theorem}  (\cite{Bray_Penrose}) 
{\em 
Let $(M^3, g)$ be a complete, asymptotically flat $3$-manifold with boundary with
nonnegative scalar curvature. Suppose its boundary $\partial M$ consists of horizons 
(that is $\partial M$ has zero mean curvature). Let $G(x)$ be a function on $M^3$ 
which satisfies
\be \nonumber
\left\{
\begin{array}{rcl}
\lim_{x \rightarrow \infty} G & = & 1 \\
\triangle G & = & 0 \\
G|_{\partial M} & = & 0  .
\end{array}
\right.
\ee
Then
\be 
m \geq  C,
\ee
where $m$ is the total mass of $(M^3, g)$ and $C$ is the constant
in the asymptotic expansion 
\be \nonumber 
G(x) = 1 - \frac{C}{|x|} + O \left( \frac{1}{|x|^2} \right) \ \ \mathrm{at}
\ \infty.
\ee
Furthermore, equality holds if and only if the manifold 
$(M^3, g)$ is isometric to a spatial Schwarzschild manifold outside 
its horizon.
}

\vh

This theorem was established to prove the Riemannian Penrose Inequality 
in \cite{Bray_Penrose}. It was later applied in \cite{Miao_hbdry} for a generalization 
of Bunting and Masood's rigidity theorem \cite{Bunting_Masood} on static vacuum 
spacetime with black hole boundary. Considering these applications, 
it is of interest to know whether a simiar theorem holds on an asymptotically 
flat $3$-manifold with a boundary that does not necessarily have zero mean curvature. 
A corollary to our main result, Theorem \ref{mainthm} stated in a moment,  
is the following:

\begin{cor} \label{maincor}
Let $(M^3, g)$ be a complete, asymptotically flat $3$-manifold with nonnegative
scalar curvature with a connected smooth boundary.  Assume $(M^3, g)$ is 
diffeomorphic to $\mathbb{R}^3 \setminus \Omega$,  where $\Omega$ is a 
bounded domain. Let $\G(x)$ be a function on $(M^3, g)$ which satisfies
\be \nonumber
\left\{
\begin{array}{rcl}
\lim_{x \rightarrow \infty} \G & = & 1 \\
\triangle \G & = & 0 \\
\G|_{\partial M} & = & \sqrt{\frac{1}{16 \pi}\int_{\partial M} H^2 d\mu }  ,
\end{array}
\right.
\ee
where $H$ is the mean curvature of $\partial M$ and $d \mu$ is the 
induced surface measure. If $\partial M$ has nonnegative Hawking mass 
(that is $\int_{\partial M} H^2 d \mu \leq 16 \pi$), then
\be 
m \geq  \C,
\ee
where $m$ is the total mass of $(M^3, g)$ and $\C$ is the constant
in the asymptotic expansion 
\be \nonumber 
\G(x) = 1 - \frac{\C}{|x|} + O \left( \frac{1}{|x|^2} \right) 
\ \ \mathrm{at}
\ \infty.
\ee
Furthermore, equality holds if and only if the manifold 
$(M^3, g)$ is isometric to a spatial Schwarzschild manifold 
$$
(M_{r_0}, g^S_m) = \left([r_0, \infty) \times S^2, \frac{1}{1 - \frac{2m}{r}} dr^2 
+ d\sigma^2 \right),
$$
where $r_0$ is some positive constant satisfying $r_0 \geq 2m$ 
and $d\sigma^2$ is the standard metric on the unit spere 
$S^2 \subset \mathbb{R}^3$. 
\end{cor}

Before we state our main result, we first define the capacity of a surface. 

\begin{df}
Let $(M^3, g)$ be a complete, asymptotically flat $3$-manifold with
a nonempty boundary $\S$. The {\bf capacity} of $\S$ in $(M^3, g)$, 
denoted by $C_M(\S, g)$, is defined to be 
\be
C_M(\S, g) = \inf \left\{ \frac{1}{4 \pi} \int_{M^3} 
| \nabla \phi |^2 dg \right\} ,
\ee
where the infimum is taken over all locally Lipschitz $\phi(x)$ which go
to $1$ at $\infty$ and equal $0$ on $\S$.
\end{df}

When $(M^3, g)$ is the complement of a smooth bounded domain $\Omega$
in the Euclidean space $(\R^3, g_0)$, $C_{M}(\partial \Omega, g_0)$ is simply 
the usual electrostatic capacity of $\partial \Omega$ \cite{Polya_Szego}. 
In this case, we write $C_{M}(\partial \Omega, g_0)$ as $C(\partial \Omega)$.
Our main theorem is the following:

\begin{thm} \label{mainthm}
Let $(M^3, g)$ be a complete, asymptotically flat $3$-manifold with nonnegative
scalar curvature with a connected smooth boundary. 
Assume $(M^3, g)$ is diffeomorphic to $\mathbb{R}^3 \setminus \Omega$, 
where $\Omega$ is a bounded domain. Then
\be \label{estofc}
C_M(\partial M, g) \leq 
\sqrt{\frac{|\partial M|}{16 \pi}} \left( 1 + \sqrt{\frac{1}{16 \pi}
\int_{\partial M} H^2 d\mu} \right),
\ee
where $| \partial M |$ and $H$ are the area and the mean 
curvature of $\partial M$. 
Furthermore, equality holds if and only 
$(M^3, g)$ is isometric to a spatial Schwarzschild manifold 
$$
(M_{r_0}, g^S_{m}) = 
\left([r_0, \infty) \times S^2, \frac{1}{1 - \frac{2m}{r}} dr^2 + 
d\sigma^2 \right),
$$
where $r_0$ is some positive constant satisfying $r_0 \geq 2 m$ and 
$d\sigma^2$ is the standard metric on the unit sphere $S^2 \subset
\mathbb{R}^3$.
\end{thm}

As an immediate corollary, we have a new estimate of the 
capacity of a surface in $(\R^3, g_0)$. 
 
\begin{cor} 
Let $\Omega \subset (\R^3, g_0)$ be a bounded domain with a 
connected smooth boundary. Then 
\be \label{flatest}
C(\partial \Omega) \leq \sqrt{\frac{|\partial \Omega|}{16 \pi}} 
\left( 1 + \sqrt{\frac{1}{16 \pi} \int_{\partial \Omega} H^2 d\mu} \right).
\ee
Furthermore, equality holds if and only if $\Omega$ is a round ball.
\end{cor}

We note that it is interesting to compare (\ref{flatest}) 
with the classical isoperimetric inequality in $\R^3$:
\be \label{classiciso}
\left(  \frac{3 V}{4 \pi} \right)^\frac{1}{3}
\leq \sqrt{ \frac{ | \partial \Omega |}{4 \pi}} 
=  \sqrt{\frac{|\partial \Omega|}{16 \pi}} \left( 1 + 1  \right) ,
\ee
where $V$ is the volume of $\Omega \subset \R^3$.
It is known, among all domains $\Omega$ with a fixed 
amount volume $V > 0$, $C(\partial \Omega)$ is minimized by a round ball 
\cite{Polya_Szego}, i.e.
\be
\left(  \frac{3 V}{4 \pi} \right)^\frac{1}{3} \leq  C(\partial \Omega ) .  
\ee
On the other hand,  the Willmore functional 
$\int_{\partial \Omega} H^2 d\mu$ satisfies
$$ \int_{\partial \Omega} H^2 d\mu \geq 16 \pi .$$
Hence, (\ref{flatest}) is analogous to the classical isoperimetric 
inequality (\ref{classiciso}), but where both sides of the inequality are increased.

We outline the idea of the proofs of Corollary \ref{maincor} and Theorem 
\ref{mainthm}. For both results, we apply the technique of weak inverse 
mean curvature flow as developed by Huisken and Ilmanen in 
\cite{IMF}. The topological assumptions on $M^3$ 
ensures that the Hawking mass of the flowing surface $\S_t$ 
is monotone nondecreasing for positive $t$. 
A key step in the proof of Theorem \ref{mainthm} is to use the flow to construct a special 
test function that gives the estimate (\ref{estofc}). Such a construction was 
first used by Bray and Neves in \cite{Bray_Neves}. 
It is a convenient feature of Corollary \ref{maincor} and Theorem \ref{mainthm} that,
though they are proved by applying the weak inverse mean curvature flow technique, 
they hold {\em without} assuming the boundary surface
is {\em outer minimizing} (see \cite{Bray_Penrose} \cite{IMF}).

This paper is organized as follows. In Section \ref{literature}, we recall some 
classical results on capacity of convex surfaces in $\R^3$ from the 
work of P\'olya and Szeg\"o in \cite{Polya_Szego} to illustrate that the weak 
inverse mean curvature flow technique fits naturally with the classical method. 
In Section \ref{ucap}, we apply the theory of weak inverse mean curvature 
flow to prove a general theorem on $C_M(\S, g)$, which 
includes Thereom \ref{mainthm} as a special case. 
In Section \ref{application}, we relate our estimate on $C_M(\S,g)$ to estimates 
of the Hawking mass and the total mass, and prove Corollary \ref{maincor}. 
In Section \ref{discussion}, we give an application of Corollary \ref{maincor} to 
the study of static metrics in general relativity.

\section{Capacity of convex surfaces in $\R^3$} \label{literature}
We first give an account of  some classical methods and results from $\cite{Polya_Szego}$
in estimating $C(\S)$ for convex surfaces $\S$ in $\R^3$. 

Let $\S$ be a closed, connected  $C^2$ surface bounding some domain $\Omega$
in $\R^3$. One basic idea in estimating $C(\S)$ is to minimize $\int_{\R^3 \setminus
\Omega} | \nabla v |^2 dg_0$ over functions 
$v$ which have {\em given level surfaces} $\{ \S_t \}$. These level surfaces form 
a one-parameter family. Therefore, after the selection of $\{\S_t\}$, $v$ becomes a 
function of one variable and the infimum of $\int_{\R^3 \setminus \Omega}  
| \nabla v|^2 d g_0$ over all such $v$ 
can be easily evalutated. Precisely, we fix a function $\psi$, defined on 
$\R^3 \setminus \Omega$, which
satisfies
$$ \psi \geq 0, \ \ \ \psi|_\S = 0 , \ \ \ \lim_{x \rightarrow \infty} \psi = \infty .$$
Let  
$
\{ \S_t = \psi^{-1}(t) \ | \  0  \leq t < \infty \} 
$
be the family of level surfaces of $\psi$.
For any other function $v $ having the same level surfaces $\{ \S_t \}$, 
$v$ must have the form $ v(x) = f( \psi(x) )$ for some single variable function $f(t)$,
which satisfies $f(0) = 0$ and  $f(\infty) = 1$. By the co-area formula, we have
\beq
\int_{\R^3 \setminus \Omega} | \nabla v |^2 dg_0 & = & 
\int_0^\infty \left(   \int_{ \S_t  } f^\prime(t)^2 | \nabla \psi |  d \mu 
\right) dt  \nonumber \\
& = & \int_0^\infty f^\prime(t)^2   \left(   \int_{ \S_t } | \nabla \psi |  d \mu 
\right) dt .
\eeq
Define \label{defTt}
\be
T(t) = \frac{1}{4 \pi} \int_{ \S_t } | \nabla \psi | d \mu,
\ee
which is determined solely by the level surfaces $\{\S_t\}$. Then
\be 
C(\S) \leq \int_0^\infty f^\prime(t)^2 T(t) dt.  
\ee  
Applying the fundamental theorem of calculus and the H\"older inequality, 
\beq
1 & = & \left( \int_0^\infty f^\prime(t) dt \right)^2 =
\left( \int^\infty_0 f^\prime(t) T(t)^\frac{1}{2} T(t)^{-\frac{1}{2}} dt \right)^2 \nonumber \\
& \leq & \left( \int^\infty_0 f^\prime(t)^2 T(t) dt \right)
\left( \int_0^\infty  T(t)^{-1} dt  \right) .
\eeq 
Thus 
\be
 \left( \int_0^\infty  T(t)^{-1} dt  \right)^{-1}  \leq \int_0^\infty f^\prime(t)^2 T(t) dt
\ee
for all $f(t)$ with equality if and only if 
\be
f(t) = \Lambda \int_0^t \frac{1}{T(s)} ds,
\ee
where $\Lambda =  \left( \int_0^\infty  T(t)^{-1} dt  \right)^{-1}$. Choosing such a $f(t)$,  
we show that
\be
C(\S) \leq  \left( \int_0^\infty  T(t)^{-1} dt  \right)^{-1} .
\ee

Now suppose $\S$ is a {\em convex} surface in $\R^3$. A natural choice for 
$\{ \S_t \}$ is the family of level surfaces of the distance function to $\S$.  We 
let $\psi(x) = dist(x, \S)$ and define
\be
\S_t = \{ x \ | \ dist(x, \S) = t \}.
\ee
Then $ | \nabla \psi | = 1$ everywhere and $ T(t) = \frac{ | \S_t | }{4 \pi} $,
where $| \S_t |$ is the area of $\S_t$, given by 
\be
| \S_t | = | \S | + \left( \int_{\S} H d \mu \right) t + 4 \pi t^2 .
\ee
This leads to the following theorem of Szeg\"o.

\vh

\noindent {\bf Theorem} (\cite{Polya_Szego})
{\em If $\S$ is a convex surface in $(\R^3, g_0)$, then
\be
C(\S) \leq \frac{M}{ 4 \pi} \frac{2 \ep}{\log{ \frac{1 + \ep}{1 - \ep}}}, 
\ee
where 
\be
M = \frac{1}{2} \int_\S  H d \mu \ \mathrm{and} \  
\ep^2 = 1 - \frac{4 \pi |\S|}{ M^2}.
\ee
}
 
\vh

\section{Capacity of surfaces in asymptotically flat manifolds} \label{ucap}
Obviously, most of the calculations in Section 2 work on any non-compact 
Riemannian manifold. The key step is to make a good choice of $\{ \S_t \}$
so that the corresponding $T(t)$  can be efficiently estimated. 
In this section, we consider {\em asymptotically flat $3$-manifolds}, on which 
the theory of weak inverse mean curvature flow developed by Huisken and Ilmanen
\cite{IMF} gives a nearly canonical foliation.

\begin{df}
A Riemannian $3$-manifold $(M^3, g)$ is said to be {\bf asymptotically flat}
if there is a compact set $K \subset M$ such that $M \setminus K$ is 
diffeomorphic to $\R^3$ minus a compact set  and in the 
coordinate chart defined by this diffeomorphism, 
$$
g = \sum_{i,j} g_{ij}(x) dx^i dx^j,
$$
where 
$$ g_{ij} = \delta_{ij} + O(|x|^{-1}), \ \partial_k g_{ij} = O(|x|^{-2})
\ \ \partial_l \partial_{k} g_{ij} = O(|x|^{-3}) .$$
\end{df}

\vh

\noindent {\bf Theorem} (\cite{IMF})
{\em
Let $(M^3, g)$ be a complete, connected asymptotically flat $3$-manifold 
with a $C^{1,1}$ boundary 
$\S$.  There exists a proper, locally Lipschitz 
function $\phi \geq 0$ on $M$, called the solution to the weak inverse mean curvature
flow with initial condtion $\S$, which satisfies the following properties:
\begin{enumerate}

\item $\phi|_\S = 0$, $\lim_{x \rightarrow \infty} \phi = \infty$. 
For $t > 0$, $\S_t = \partial \{ \phi \geq t\} $ and $\S^\prime_t =
\partial \{ \phi > t \}$
define an increasing family of $C^{1, \alpha}$ surfaces.

\item The surfaces $\S_t$ ($\Sp_t$) minimize (strictly minimize) 
area among surfaces homologous to $\S_t$ in the region $\{ \phi \geq t \}$. 
The surface $ \Sp = \partial \{ \phi >0 \}$
strictly minimizes area among surfaces homologous to $\S$ in $M$.

\item For almost all $t > 0$,  the weak mean curvature of $\S_t$ is defined
and equals $| \nabla \phi |$, which is positive for almost all $ x \in  \S_t $.

\item  For each $t > 0$, $|\S_t |  = e^t | \S^\prime| $, and
$| \S_t | = e^t | \S |$ if $\S$ is outer minimizing (that is $\S$ minimizes area among all
surfaces homolgous to $\S$ in $M$).

\item  If $(M^3, g)$ has nonnegative scalar curvature and $\chi(\S_t) \leq 2$ for 
all $t > 0$, the Hawking mass 
$$ m_H(\S_t) = \sqrt{ \frac{| \S_t|}{16 \pi} } \left( 1 - \frac{1}{16 \pi} 
\int_{\S_t} H^2 d \mu \right) $$ 
is monotone nondecreasing for $t > 0$ and 
$lim_{t \rightarrow 0+} m_H(\S_t) \geq m_H(\S^\prime)$.  
Here $\chi(S)$ is the Euler characteristic of a surface $S$. 

\end{enumerate}
}

\vh

We note that when $\phi$ is a smooth function with non-vanishing gradient, 
property $3$ is just saying that the level surfaces $\{ \S_t \}$ move at a 
speed equal to the inverse of their mean curvature.
We are now ready to prove the main result of this section.

\begin{thm} \label{gthm}
Let $(M^3, g)$ be a complete, connected asymptotically flat $3$-manifold which has a $C^2$
boundary $\S$. Suppose $(M^3, g)$ has nonnegative scalar curvature and the solution to 
the weak inverse mean curvature flow with initial condition $\S$ 
satisfies $\chi(\S_t) \leq 2$ for all $t > 0$. Then
\be \label{estofcL}
C_M(\S, g) \leq 
\sqrt{\frac{|\S|}{16 \pi}} \left( 1 + \sqrt{\frac{1}{16 \pi}
\int_\S H^2 d\mu} \right),
\ee
where $| \S |$ and $H$ denote the area and the mean 
curvature of $\S$. 
Furthermore, equality holds if and only 
$(M^3, g)$ is isometric to a spatial Schwarzschild manifold 
\be \label{schwarzschild}
(M_{r_0}, g^S_{m}) = 
\left([r_0, \infty) \times S^2, \frac{1}{1 - \frac{2m}{r}} dr^2 + 
d\sigma^2 \right),
\ee
where $r_0$ is some positive constant satisfying $r_0 \geq 2 m$ and 
$d\sigma^2$ is the standard metric on the unit sphere $S^2 \subset
\mathbb{R}^3$.

\end{thm}

\pf We estimate $\int_M | \nabla v |^2 dg$ for functions $v$ that 
have the same level surfaces $\{ \S_t \}$ as the function $\phi$, which
is the solution to the weak inverse mean curvature flow starting at $\S$.
It follows from the calculation in Section 2 that
\be \label{fT}
C_M(\S, g) \leq \inf_f \left\{ \int_0^\infty f^\prime(t)^2 T(t) dt \right\} ,   
\ee
where the infimum is taken over all $f(t)$ satisfying $f(0)=0$ and $f(\infty) = 1$
and
\be
T(t) = \frac{1}{ 4 \pi} \int_{\S_t} | \nabla \phi | d\mu .
\ee
Now, for a.e. $t>0$,
\be
\int_{\S_t} | \nabla \phi | d\mu = \int_{\S_t} H d \mu,
\ee
where $H$ is the weak mean curvature of $\S_t$. To proceed, we make use of 
the key property that $m_H(\S_t)$ is monotone nondecreasing for $t>0$
and $m_H(\Sp) \leq \lim_{s \rightarrow 0+} m_H(\S_s) $,
which are guaranteed by the assumption that 
$(M^3, g)$ has nonnegative scalar curvature and 
$\chi(\S_t) \leq 2$. Hence, for each $t > 0$, 
\be
m_H(\Sp) \leq  m_H(\S_t) = \sqrt{\frac{|\S_t|}{16 \pi}} 
\left( 1 - \frac{1}{16 \pi} \int_{\S_t} H^2 d\mu
\right) .
\ee
This, together with the H\"older inequality, implies
\be
\frac{1}{16 \pi |\S_t| } \left( \int_{\S_t} H d\mu \right)^2 
 \leq \frac{1}{16 \pi}
\int_{\S_t} H^2 d\mu \leq 1 - m_H(\S^\prime) \sqrt{ \frac{16 \pi}{|\S_t|} }.
\ee
Hence,
\beq \label{b}
4 \pi T(t) = \int_{\S_t} H d\mu  & \leq & 
\left[ 16 \pi | \S_t |   \left( 1 - m_H(\S^\prime) \sqrt{ \frac{16 \pi}{|\S_t|} }
\right) \right]^{\frac{1}{2}} \nonumber \\
& = & \left[ 16 \pi | \S^\prime | e^t  \left( 1 -  m_H(\S^\prime) \sqrt{ \frac{16 \pi}{|\S^\prime|} } 
e^{-\frac{t}{2}} 
\right) \right]^{\frac{1}{2}} , 
\eeq
by the fact $| \S_t | = e^t |\S^\prime|$.
Now write $\mpz = m_H(\Sp) $ and $ \Abz = | \Sp | $, it follows from (\ref{fT}) and (\ref{b}) that
\be \label{cc}
4 \pi C_M(\S, g)   \leq \inf_f \left\{ \int_0^\infty f^\prime(t)^2 F( \Abz, \mpz, t) dt \right\},
\ee
where $F(\Abz, \mpz, t)$ is an explicit function of $\Abz$, $\mpz$ and $t$, given by
\be
F(\Abz, \mpz, t) =  \left[ 16 \pi \Abz e^t  \left( 1 - \mpz 
\sqrt{ \frac{16 \pi}{\Abz} } e^{-\frac{t}{2}} \right) \right]^{\frac{1}{2}}.
\ee

To calculate 
$$
\inf_f \left\{ \int_0^\infty f^\prime(t)^2 F( \Abz, \mpz, t) dt \right\},
$$
we consider the $3$-dimensional spatial 
Schwarzschild metric 
$$ g^S_{\mpz} = \frac{1}{1 - \frac{2 \mpz}{r}} dr^2 + r^2 d\sigma^2 . $$
When $\mpz < 0$, $g^S_{\mpz}$ is defined on $(0, \infty) \times S^2$ 
(the metric has a  singularity  at $r=0$). 
When $\mpz \geq 0$, $g^S_{\mpz}$ is defined on $[2 \mpz, \infty) \times S^2$. 
In either case, 
$g^S_{\mpz}$ is well defined on $[\rbz, \infty) \times S^2$, where
$\rbz$ satisfies  $\Abz = 4 \pi \rbz^2$. For convenience, we let 
$M_{\rbz}$ denote $[\rbz, \infty) \times S^2$, then the spatial 
Schwarzschild manifold $(M_{\rbz}, g^S_{\mpz})$ has a boundary 
$\{r = \rbz\}$ whose area is $\Abz$. A basic fact about $(M_{\rbz}, g^S_{\mpz})$
is that the classic inverse mean curvature flow in $(M_{\rbz}, g^S_{\mpz})$ with 
initial condition $\{ r  = \rbz  \}$ is given by the family of coordinate spheres
\be
S_t = \{ r =  \rbz e^{\frac{1}{2}t} \} , 
\ee
which have constant mean curvature (depending on $t$) and  
constant Hawking mass $\mpz$. 
Therefore, the corresponding function $\phi(x)$ on $(M_{\rbz}, g^S_{\mpz})$ 
is given by 
\be 
\phi(x) = 2 \log{ \left( \frac{r}{\rbz}  \right) } .
\ee 
Next, we consider the harmonic function
$u$ on $(M_{\rbz}, g^S_{\mpz})$ that equals $0$ at $\{r = \rbz\}$ and goes to $1$ 
at $\infty$. We have two cases:

\vspace{.2cm} 

\noindent {\bf Case 1} $\mpz \neq 0$:

\vspace{.2cm}

\noindent In this case, the function 
\be
v = \sqrt{1 - \frac{2 \mpz }{r} } = 1 - \frac{\mpz}{r} + O \left(\frac{1}{r^2}\right)
\ee
is a non-constant harmonic function on $(M_{\rbz}, g^S_{\mpz})$. Hence, 
\be
u =  \frac{v - v(\rbz)}{ 1 - v(\rbz) }.
\ee
The capacity of the boundary of $(M_{\rbz}, g^S_{\mpz})$ is 
\be \label{modelcap}
C_{M_{\rbz}}(\partial M_{\rbz}, g^S_{\mpz}) = 
\frac{1}{4\pi} \int_{M_{\rbz}} | \nabla {u} |^2 d g^S_{\mpz}
= \frac{\mpz}{1 - v(\rbz)}.
\ee
Note that $u$ can be rewritten as 
\be \label{fphi}
u = f_0 \circ \phi ,
\ee
where
\be \label{fa}
f_0(t) = \frac{1}{1 - v(\rbz) } \left[ 
\sqrt{ 1 - \frac{2 \mpz }{ \rbz e^{\frac{t}{2}}}} - v(\rbz)
\right] .
\ee 
It then follows from (\ref{modelcap}), (\ref{fphi}) and the fact  that $S_t$ has constant 
mean curvature and  $m_H(S_t) = \mpz$ for all $t$  that $f_0$ achieves
$$
\inf_f \left\{ \int_0^\infty f^\prime(t)^2 F( \Abz, \mpz, t) dt \right\}
$$
and the infimum is given by
\beq
\int_{M_{\rbz}} | \nabla u |^2 & = & 
4 \pi \frac{\mpz}{1 - v(\rbz)} \nonumber \\
& = & 4 \pi \sqrt{\frac{|\Sp|}{16 \pi}} \left( 1 + \sqrt{\frac{1}{16 \pi}
\int_{\Sp} H^2 d\mu} \right) .
\eeq
Going back to (\ref{cc}), we have
\be \label{aa}
C_M( \S,  g) \leq  \sqrt{\frac{|\Sp|}{16 \pi}} \left( 1 + \sqrt{\frac{1}{16 \pi}
\int_{\Sp} H^2 d\mu} \right) .
\ee

\vspace{.2cm}

\noindent {\bf Case 2} $\mpz = 0$:

\vspace{.2cm}

\noindent In this case, our model space $(M_{\rbz}, g^S_{\mpz})$
is the Euclidean space $(\R^3, g_0)$ minus a round ball of radius 
$\rbz$ centered at the origin. Hence, 
\be
u = 1 - \frac{\rbz}{r}. 
\ee
The capacity of the boundary of $(M_{\rbz}, g^S_{\mpz})$ is 
\be
C(\partial M_{\rbz}) = \rbz . 
\ee
Defining 
\be \label{fb}
f_0 (t) = 1 - e^{- \frac{t}{2}} ,
\ee
we can rewrite  $u$ as 
\be
 u = f_0 \circ \phi . 
\ee
The same argument as in the Case 1 then implies that
\be \label{bb}
C_M( \S,  g) \leq  \rbz = \sqrt{\frac{|\Sp|}{4 \pi}} 
= \sqrt{\frac{|\Sp|}{ 16 \pi}}
\left( 1 + \sqrt{\frac{1}{16 \pi}
\int_{\Sp} H^2 d\mu} \right) ,
\ee
where the last equality holds because $\mpz = 0 $. 

\vspace{.2cm} 

Therefore, in both cases, we have proved that 
\be \label{dd}
C_M( \S,  g) \leq  
\sqrt{\frac{|\Sp|}{ 16 \pi}}
\left( 1 + \sqrt{\frac{1}{16 \pi}
\int_{\Sp} H^2 d\mu} \right) .
\ee
To replace $\S^\prime$ by $\S$,  we use the property  
that $\Sp$ strictly minimizes area among all surfaces homologous to $\S$. 
Since $\S$ is $C^2$, $\Sp$ is $C^{1,1}$ and $C^\infty$ where 
$\Sp$ does not contact $\S$. 
Moreoever, the mean curvature $H^\prime$ of 
$\Sp$ satisfies
\be 
H^\prime = 0 \ \mathrm{on} \ \Sp \setminus \S \ \ \ \mathrm{and} \ \ \
H^\prime = H \geq 0 \ \ \mathcal{H}^2 a.e. \ \mathrm{on} \  \Sp \cap \S .
\ee
In particular, we have
\be \label{sandsp}
|\Sp| \leq |\S| \ \ \mathrm{and} \ \
 \int_{\Sp} {H^\prime}^2 d \mu \leq \int_\S H^2 d\mu .
\ee
Therefore, it follows from (\ref{dd}) and (\ref{sandsp}) that 
\be \label{kk}
C_M(\S, g) \leq   \sqrt{\frac{|\S |}{16 \pi}} \left( 1 + \sqrt{\frac{1}{16 \pi}
\int_{\S} H^2 d\mu} \right) .
\ee

\vspace{.2cm}

To complete the proof of Theorem \ref{gthm}, we must consider the case of
equality. Suppose
\be \label{equality}
C_M(\S, g)  =  
\sqrt{\frac{|\S|}{16 \pi}} \left( 1 + \sqrt{\frac{1}{16 \pi}
\int_\S H^2 d\mu} \right) ,
\ee
it follows from the above proof that
\be \label{outm}
|\Sp| = |\S|, \   
 \int_{\Sp} {H^\prime}^2 d \mu =  \int_\S H^2 d\mu  , 
\ee
\be \label{hconstant}
m_H(\S_t) = m_H(\Sp), \ \forall \ t>0 ,
\ee
and
\be \label{capvalue}
C_M(\S, g) = \frac{1}{4\pi}  \int_M | \nabla ( f_0 \circ \phi) |^2 dg,
\ee
where $f_0$ is either given by (\ref{fa}) or (\ref{fb}) and 
$\phi$ is the solution to the weak inverse mean curvature flow 
in $(M^3, g)$ with initial condition $\S$. 
It follows from (\ref{outm}) and (\ref{hconstant}) that
$\S$ is outer minimizing and the Hawking mass
$m_H(\S_t)$ equals $m_H(\S)$ for every $t$.  On the other hand,
(\ref{capvalue}) implies that 
\be
u_M = f_0 \circ \phi 
\ee 
is a smooth harmonic function on $M$. 
As a result, the surfaces $\S_t$, $\S$ do not ``jump" to $\S_t^\prime$, $\S^\prime$ (as 
defined in \cite{IMF}, meaning $\S_t = \S_t^\prime, \ \S = \Sp$), for 
otherwise the set $\{ x \in M^3 \ | \ u_M(x) = f_0(t) \} $ would have non-empty interior
for some $t \geq 0$, contradicting the maximum principle for harmonic functions.
Furthermore, applying the maximum principle to the exterior region of $\S_t$ in $(M^3, g)$
and using the fact that $u_M$ is constant on $\S_t$ and $\S_t$ is at lease $C^1$, 
we conclude that $\nabla u_M$ never vanishes. 
Therefore, 
\be
\phi = f^{-1}_0 \circ u_M
\ee
is a smooth function on $M$ with non-vanishing gradient. Hence, $\{\S_t\}$ evolve 
smoothly at a speed equal to the inverse of their mean curvature.
The fact $m_H(\S_t) = m_H(\S)$ for all $t>0$ then readily implies
that $(M^3, g)$ is isometric to a spatial Schwarzschild manifold
\be
(M_{r_0}, g^S_{m}) = 
\left([r_0, \infty) \times S^2, \frac{1}{1 - \frac{2m}{r}} dr^2 + 
d\sigma^2 \right) 
\ee
with $r_0 \geq 2m$ (see page 423 in \cite{IMF} for detail).
Theorem \ref{gthm} is proved. \stop

\vspace{.3cm}

Next, we give a topological condition 
of $M^3$ that is sufficient to guarantee
the assumption $\chi(\S_t) \leq 2$ in Theorem \ref{gthm}.

\begin{prop} \label{topp}
Let $(M^3, g)$ be a complete, connected asymptotically flat $3$-manifold with a nonempty
boundary $\S$. If $H_2(M,\S) = 0$ and $\S$ is connected, then 
$\{\S_t\}_{t>0}$ remains connected. In particular, $\chi(\S_t) \leq 2$ for all $t$. 
\end{prop}

\pf  Under the assumption that $\S$ is connected, Huisken and Ilmanen proved that 
the sets $\{ \phi < t \}$ and $\{ \phi > t \}$ are connected \cite{IMF}. They also showed,
for each $t >0$, $\S_t$ can be approximated in $C^1$ by earlier surfaces $\S_s$,  
satisfying $\S_s = \S_s^\prime$. Hence, we may assume $\S_t = \{ \phi = t \}$.  

Let $\S_1$ be one component of $\S_t$. Since $H_2(M,\S) =0$
and $\S$ is connected, there is a bounded region $D$ in $M$ such that 
either $\partial D = \S_1$ or
$\partial D = \S_1 \cup \S$. As the set  $\{ \phi > t \}$ is connected and contains 
$\infty$, we have $\phi \leq t$ on $D$. It then follows that $\phi < t$ on
$D \setminus \S_1$. Hence $D \setminus \S_1$ is a component of the set $\{ \phi < t \}$.
Since $\{ \phi < t \}$ is connected, we must have $D \setminus \S_1 = \{ \phi < t \}$.
Therefore, $\S_1$ is the only component of $\S_t$ (and  
$\partial D = \S_1 \cup \S$).  We conclude that $\S_t$ is connected.  \stop

\vspace{.3cm} 

Theorem \ref{mainthm} now follows directly from Theorem \ref{gthm}, Proposition
\ref{topp} and the fact $H_2(\R^3 \setminus \Omega, \partial \Omega) = 0$ . 

\section{Estimate of the total mass} \label{application}
We prove Corollary \ref{maincor} in this section. 
First, we point out that Theorem \ref{gthm} translates directly into a statement 
about the capacity of $\S$ and the Hawking mass of $\S$. 

\begin{thm} \label{Hmass}
Let $(M^3, g)$ be a complete, connected asymptotically flat $3$-manifold with a 
nonempty boundary $\S$. Suppose $(M^3, g)$ satisfies all the assumptions
in Theorem \ref{gthm}. Then
\be \label{Hmasscap}
 | m_H(\S) | \geq | 1 - \alpha | C_M(\S, g) ,
\ee
where $\alpha$ is a constant defined by 
\be
\alpha = \sqrt{ \frac{1}{16\pi} \int_\S H^2 d \mu }.
\ee
Furthermore,  in the case $ \alpha  \neq 1$, equality holds if and only if   
$(M^3, g)$ is isometric to a spatial Schwarzschild manifold 
\be \label{isometry}
(M_{r_0}, g^S_{m}) = 
\left([r_0, \infty) \times S^2, \frac{1}{1 - \frac{2m}{r}} dr^2 + 
d\sigma^2 \right),
\ee
where $r_0$ is some positive constant satisfying $r_0 \geq 2m$ and 
$d\sigma^2$ is the standard metric on the unit sphere $S^2 \subset
\mathbb{R}^3$.
\end{thm}

\pf It follows directly from Theorem \ref{gthm} and the definition of the
Hawking mass. \stop

\vspace{.3cm}

Next, we recall the definition of the total mass of an asymptotically flat $3$-manifold 
(see \cite{Bray_Penrose}, \cite{IMF}).

\begin{df}
Let $(M^3, g)$ be an asymptotically flat $3$-manifold. The total mass of $(M^3, g)$
is defined as the limit
\be
m = \frac{1}{16 \pi} \lim_{r \rightarrow \infty}
\int_{S_r} \sum_{ij} ( \partial_i g_{ij} - \partial_j g_{ii} ) \nu^j d \sigma ,
\ee
where $S_r$ is the coordinate sphere $\{ |x| = r\}$, $\nu$ is the coordinate unit normal
to $S_r$  and $ d \sigma$ is the area element of $S_r$ in the coordinate chart.
\end{df}

An important feature of the weak inverse mean curvature flow, 
proved by Huisken and Ilmanen in \cite{IMF}, is 

\begin{prop} \label{IMFlimit}
Let $(M^3, g)$ be an asymptotically flat $3$-manifold and $\phi$
be a solution to the weak inverse mean curvature flow, then
\be 
m \geq \lim_{t \rightarrow \infty} m_H(\S_t) ,
\ee
where $m$ is the total mass of $(M^3, g)$ and $\S_t = \partial \{ \phi \geq t\}$.
\end{prop}

The next theorem shows that the total mass $m$ is bounded from below 
by the same quantity $( 1 - \alpha) C_M(\S, g)$ as in (\ref{Hmasscap}). 
A convenient feature of the theorem is that it does not require $\S$ to 
be outer minimizing. 

\begin{thm}  \label{totalmass} 
Let $(M^3, g)$ be a complete, connected asymptotically flat $3$-manifold 
with a nonempty boundary $\S$.   
Suppose $(M^3, g)$ satisfies all the assumptions in Theorem \ref{gthm}
and $\S$ has nonnegative Hawking mass. 
Let $\G(x)$ be a function on $(M^3, g)$ which 
satisfies
\be \nonumber
\left\{
\begin{array}{rcl}
\lim_{x \rightarrow \infty} \G & = & 1 \\
\triangle \G & = & 0 \\
\G|_{\S} & = & \alpha ,
\end{array}
\right.
\ee
where 
\be
\alpha =  \sqrt{\frac{1}{16 \pi}\int_\S H^2 d\mu } . 
\ee
Then
\be 
m \geq  \C,
\ee
where $m$ is the total mass of $(M^3, g)$ and $\C$ is the constant
in the asymptotic expansion 
\be \nonumber 
\G(x) = 1 - \frac{\C}{|x|} + O \left( \frac{1}{|x|^2} \right) 
\ \ \mathrm{at}
\ \infty.
\ee
Furthermore, the equlity holds if and only if the manifold 
$(M^3, g)$ is isometric to a spatial Schwarzschild manifold 
$$
(M_{r_0}, g^S_m) = \left([r_0, \infty) \times S^2, \frac{1}{1 - \frac{2m}{r}} dr^2 
+ d\sigma^2 \right),
$$
where $r_0$ is some positive constant satisfying $r_0 \geq 2m$ and $d\sigma^2$ is the standard 
metric on the unit spere $S^2 \subset \mathbb{R}^3$. 
\end{thm}

\pf 
We only need to consider the case $m_H(\S)>0$ (that is $\alpha < 1$), 
as the case $m_H(\S)=0$ is essentially 
the proof of the Positive Mass Theorem via the inverse mean curvature flow \cite{IMF}.

We use the same notations as in the proof of Theorem \ref{gthm}. 
Applying Proposition \ref{IMFlimit} and using the monotonicity of the 
Hawking mass, we have 
\be \label{mandsp}
m \geq \lim_{t \rightarrow \infty} m_H(\S_t)
\geq \lim_{t \rightarrow 0+} m_H(\S_t) \geq  m_H(\Sp) .
\ee
The proof of Theorem \ref{gthm}  shows
\be \label{YY}
C_M(\Sp, g) \leq \sqrt{ \frac{ | \Sp|}{16 \pi} } \left( 1+ 
\sqrt{ \frac{1}{16 \pi}  \int_{\Sp} H^2 d \mu } \right)  .
\ee
Let $ \alpha^\prime = \sqrt{ \frac{1}{16\pi} \int_{\S^\prime} H^2 d \mu }$, 
by (\ref{sandsp}) we have
\be \label{alal}
\alpha^\prime \leq \alpha  .
\ee
Thus, $\alpha^\prime < 1$. Hence, (\ref{YY}) is equivalent to 
\be \label{mm}
m_H (\Sp) \geq ( 1 - \alpha^\prime) C_M(\S^\prime, g) . 
\ee
On the other hand, it follows from the maximum principle that
\be \label{capcap}
\ C_M(\S^\prime, g) \geq C_M(\S, g) 
\ee
with equality if and only if $\S = \S^\prime$. Therefore, it follows from
(\ref{mandsp}), (\ref{mm}), (\ref{alal}) and (\ref{capcap}) that
\be
m \geq ( 1 - \alpha ) C_M(\S, g) 
\ee
with equality if and only if $\S = \S^\prime$ and 
\be \label {ZZ}
C_M(\S, g) = \sqrt{ \frac{ | \S |}{16 \pi} } \left( 1+ 
\sqrt{ \frac{1}{16 \pi}  \int_{\S} H^2 d \mu } \right)  .
\ee
By Theorem \ref{gthm}, (\ref{ZZ}) holds if and only if 
$(M^3, g)$ is isometric to a spatial Schwarzschild manifold $(M_{r_0}, g^S_m)$.
As $ \C = ( 1 - \alpha ) C_M(\S, g)$, Theorem \ref{totalmass} is proved. \stop

\vh

Corollary \ref{maincor} follows directly from Theorem \ref{totalmass},
Proposition \ref{topp} and the fact $H_2(\R^3 \setminus \Omega, \partial \Omega) = 0 $.

\section{Application to static metrics} \label{discussion} 
In this section, we give a simple application of Corollary \ref{maincor} to the study 
of static metrics in general relativity. 

We recall that a $3$-dimensional asymptotically flat manifold $(M^3, g)$ is called 
{\bf static} \cite{Miao_static}
if there is a positive function $N$, called the static potential function of $(M^3, g)$, 
satisfying $N \rightarrow 1$ at $\infty$ and
\be \label{staticeq}
\left\{
\begin{array}{ccl}
N Ric(g) & = & D^2 N  \\
\triangle N & = & 0 ,
\end{array}
\right.
\ee
where $D^2 N$ is the Hessian of $N$ and $Ric(g)$ is the Ricci curvature of $g$. 
It can be easily checked that $(M^3, g)$ and $N$ satisfy (\ref{staticeq}) if and only 
if the asymptotically flat spacetime metric $\bar{g} = - N^2 dt^2 + g $ solves 
the Vaccum Einstein Equation on $ M \times \R$. In particular, (\ref{staticeq}) 
implies that $(M^3, g)$ has zero scalar curvature. 

A fundamental result in the study of static, asymptotically flat manifolds 
with boundary is the following black hole uniqueness theorem, 
proved by Bunting and Masood-ul-Alam \cite{Bunting_Masood}. 

\vh

\noindent {\bf Theorem}  (\cite{Bunting_Masood}) 
{\em 
Let $(M^3, g)$ be a static, asymptotically flat manifold with a nonempty smooth boundary. 
Let $N$ be the static potential function of $(M^3, g)$. If $N$ satisfies
\be
N |_{\partial M} = 0  ,
\ee
then $(M^3, g)$ is isometric to a spatial Schwarzschild manifold outside 
its horizon.}

\vh 

The following theorem is a direct application of Corollary \ref{maincor} and the 
maximum principle. 

\begin{thm} \label{gofBM}
Let $(M^3, g)$ be a static, asymptotically flat manifold with a connected smooth
boundary $\S$. Assume that 
$(M^3, g)$ is diffeomorphic to $\R^3 \setminus \Omega$, where $\Omega$ is 
a bounded domain. 
 Let $N$ be the static
potential function of $(M^3, g)$. If $\S$ has nonnegative Hawking mass,
then 
\be \label{estofN}
 \min_{\S} N^2 \leq  \frac{1}{16 \pi} \int_{\S} H^2 d \mu . 
\ee
Furthermore, equality holds if and only if 
$(M^3, g)$ is isometric to a spatial Schwarzschild manifold 
$$
(M_{r_0}, g^S_m) = \left([r_0, \infty) \times S^2, \frac{1}{1 - \frac{2m}{r}} dr^2 
+ d\sigma^2 \right),
$$
where $r_0$ is some positive constant satisfying $r_0 \geq 2m$ 
and $d\sigma^2$ is the standard metric on the unit spere $S^2 \subset \mathbb{R}^3$. 
\end{thm}

\noindent {\bf Remark}:
If $N|_\S = 0$, the static metric system (\ref{staticeq}) implies 
that $\S$ is totally geodesic \cite{Bunting_Masood}.
Hence the equality in (\ref{estofN}) holds automatically. Thus,
Theorem \ref{gofBM} can be viewed as a partial generalization of 
Bunting and Masood-ul-Alam's theorem.

\vh 

\pf Let $\alpha = \sqrt{\frac{1}{16\pi} \int_\S H^2 d \mu}$
and let $\G$ be the function on $M$ defined by
\be \nonumber
\left\{
\begin{array}{rcl}
\lim_{x \rightarrow \infty} \G & = & 1 \\
\triangle \G & = & 0 \\
\G|_{\partial M} & = & \alpha .
\end{array}
\right.
\ee
Consider the asymptotic expansions of $\G$ and $N$, 
$$ \G  = 1 - \frac{\C}{|x|} + {O} \left( \frac{1}{|x|^2} \right), \ as \ x \rightarrow \infty $$
$$ N  = 1 - \frac{A}{|x|} + {O} \left( \frac{1}{|x|^2} \right), \ as \ x \rightarrow \infty  . $$ 
Suppose $\min_{\S} N \geq \alpha$, the strong maximum principle then implies that 
\be \label{AC}
A \leq \C 
\ee
with equality if and only if $N = \G$. 
By analyzing the static metric system (\ref{staticeq}), 
Bunting and Masood-ul-Alm in \cite{Bunting_Masood} showed that
\be \label{AM}
A = m ,
\ee
where $m$ is the total mass of $(M^3, g)$. 
 On the other hand, as $(M^3, g)$ has zero scalar curvature and 
 $m_H(\S) \geq 0$, Corollary \ref{maincor} shows that
\be \label{MC}
 m \geq \C. 
\ee
Therefore, it follows from (\ref{AC}), (\ref{AM}) and (\ref{MC}) that
$$ m = \C  $$
and $\ N = \G $. 
In particular, $(M^3, g)$ is isometric to $(M_{r_0}, g^S_m)$ by the rigidity part 
in Corollary \ref{maincor}. 
Theorem \ref{gofBM} is proved. \stop

\bibliographystyle{plain}
\bibliography{Capacity_GR}

\end{document}